\documentclass[11pt,a4paper]{article}

\usepackage[T1]{fontenc}
\usepackage[utf8]{inputenc}
\usepackage[english]{babel}
\usepackage{amsthm,amsmath,amssymb,amsfonts,mathtools,bbm}
\usepackage{booktabs,array,longtable}
\usepackage[shortlabels]{enumitem}
\usepackage{aliascnt}
\usepackage{geometry}
\usepackage[hidelinks]{hyperref}
\usepackage[capitalise]{cleveref}
\usepackage[babel]{microtype}
\usepackage{lineno}
\usepackage{xcolor}

\newtheorem{theorem}{Theorem}[section]

\newaliascnt{proposition}{theorem}

\aliascntresetthe{proposition}

\newaliascnt{lemma}{theorem}

\aliascntresetthe{lemma}

\newaliascnt{corollary}{theorem}
\newtheorem{corollary}[corollary]{Corollary}
\aliascntresetthe{corollary}

\newaliascnt{claim}{theorem}

\aliascntresetthe{claim}

\newaliascnt{conjecture}{theorem}

\aliascntresetthe{conjecture}

\newaliascnt{problem}{theorem}
\newtheorem{problem}[problem]{Problem}
\aliascntresetthe{problem}

\theoremstyle{definition}
\newaliascnt{definition}{theorem}

\aliascntresetthe{definition}

\newaliascnt{example}{theorem}

\aliascntresetthe{example}

\theoremstyle{remark}
\newtheorem*{remark}{Remark}

\DeclareMathOperator{\VC}{VC}
\DeclareMathOperator{\Spec}{Spec}
\newcommand{\pdim}{\operatorname{dim}_{\pi}}

\DeclareMathOperator{\supp}{supp}

\newcommand{\F}{\mathbb F}
\newcommand{\Ftwo}{\mathbb F_2}
\newcommand{\R}{\mathbb R}
\newcommand{\C}{\mathbb C}
\newcommand{\one}{\mathbf 1}
\newcommand{\cube}{\{0,1\}}

\crefname{theorem}{Theorem}{Theorems}
\crefname{corollary}{Corollary}{Corollaries}
\crefname{lemma}{Lemma}{Lemmas}
\crefname{proposition}{Proposition}{Propositions}
\crefname{problem}{Problem}{Problems}
\crefname{definition}{Definition}{Definitions}
\crefname{example}{Example}{Examples}

\title{VC-Dimension vs Degree:\\ An uncertainty principle for Boolean functions}

\author{Fan Chang\thanks{School of Statistics and Data Science, Nankai University, Tianjin, China and Extremal Combinatorics and Probability Group (ECOPRO), Institute for Basic Science (IBS), Daejeon, South Korea. E-mail: 1120230060@mail.nankai.edu.cn.  Supported by the NSFC under grant 124B2019 and the Institute for Basic Science (IBS-R029-C4).} \and Yijia Fang\thanks{Department of Mathematics, National University of Singapore, Singapore. Email: fangyijia@u.nus.edu.}}
\date{}

\begin{document}

\maketitle

\begin{abstract}
We prove a support--shattering uncertainty principle for functions on the Boolean cube. Let $\mathbb{F}$ be any field and let $f:\{0,1\}^n\to\mathbb{F}$ be nonzero. If $x^S$ is a maximum-degree monomial in the multilinear representation of $f$, then $\mathrm{supp}(f)$ shatters $S^c$. Consequently,
\[
\mathrm{VC}(\mathrm{supp}(f))+\deg_{\mathbb{F}}(f)\ge n .
\]
For Boolean-valued functions this yields both the real-degree and algebraic-degree forms of the inequality.

We derive several consequences. The real Fourier support of a nonzero Boolean function satisfies
\[
\mathrm{VC}(\mathrm{Spec}(f))\ge \deg_{\mathbb{F}_2}(f),
\qquad
\mathrm{VC}(\mathrm{supp}(f))+\mathrm{VC}(\mathrm{Spec}(f))\ge n .
\]
The same principle gives a low-degree-obstruction proof of the Sauer--Shelah lemma and a product-space analogue based on the Efron--Stein decomposition, recovering the Karpovsky--Milman multivalued Sauer lemma. We also obtain polynomial-method applications, including an arbitrary-field Sziklai--Weiner lower bound and a shattering theorem for null designs, and discuss sharpness and equality examples.

\vspace{8pt}

\vspace{8pt}
\end{abstract}

\section{Introduction}

Functions defined on the Boolean hypercube $\{0,1\}^n$ are fundamental objects in Combinatorics and Theoretical Computer Science. It is well known that every such function $f:\{0,1\}^n\to\mathbb{R}$ can be represented as a linear combination,
$$
f=\sum\limits_{S\subset[n]}\hat{f}(S)\cdot \chi_S,
$$
of the $2^n$ functions $\{\chi_S\}_{S\subset[n]}$ defined by $\chi_S(x)=(-1)^{\sum_{i\in S}x_i}$. This representation is known as the \emph{Fourier expansion} of the function $f$, and the real numbers $\hat{f}(S)$ are known as its Fourier coefficients.  Considering the Fourier spectral properties of the indicator functions of set families under specific combinatorial constraints naturally establishes a bridge between the analysis of Boolean functions and Extremal Combinatorics. Tools from the analysis of Boolean functions have been successfully applied to Extremal Combinatorics in a number of works (e.g.,\cite{EFF2012,EFP2011,EKL2024,EL2022,KLLM2025,KL2021,KLM2024}).

Motivated by the study of low-influence Boolean functions, we investigate what structural properties the Fourier spectrum exhibits for Boolean functions with bounded VC-dimension. The \emph{Vapnik--Chervonenkis dimension}, or \emph{VC-dimension} for short, of a set family $\mathcal{F}\subseteq2^{[n]}$ is the maximum $d$ for which there exists a subset $S\subset[n]$ of size $d$ such that for all $T\subseteq S$, there is some $A\in\mathcal{F}$ with $T=A\cap S$. We denote the VC-dimension of a set family $\mathcal{F}$ by ${\rm VC}(\mathcal{F})$. The VC-dimension serves as a fundamental measure for quantifying the \emph{combinatorial complexity} of set systems and has found many applications in Extremal Combinatorics (e.g., \cite{BBDFP2024,CXYZ2025,FPS2019,FPS2021,FPS2023,JP2024}). 
We identify $x\in\cube^n$ with the set $\{i:x_i=1\}\subseteq[n]$ and define
\[
\supp(f):=\{x\in\cube^n:f(x)\ne0\},
\qquad
\VC(f):=\VC(\supp(f)).
\]

The paper is motivated by the following question.

\begin{problem}\label{prob:central}
Let $f:\cube^n\to\cube$ be a Boolean function with $\VC(f)=d$.  What structural information can one force on the polynomial or Fourier spectrum of $f$?
\end{problem}

Let $\F$ be a field. Every function $f:\cube^n\to\F$ has a unique multilinear representation
\begin{equation}\label{eq:multilinear}
f(x)=\sum\limits_{S\subseteq [n]}c_Sx^S,
\end{equation}
where $x^S=\prod_{i\in S}x_i$ as usual, and $c_S\in\F$. We refer to~\eqref{eq:multilinear} as the \emph{$\mathbb{F}$-polynomial representation of $f$}, and set $\deg_{\F}(f):=\max\{|S|:c_S\ne0\}$ for it \emph{$\mathbb{F}$-degree/algebraic degree}. For $\F=\R$ and Boolean-valued $f$, it agrees with the usual real, or Fourier, degree $\deg(f)$.

Our main result is the following structural theorem.

\begin{theorem}\label{thm:top-monomial}
Let $\F$ be a field and let $f:\cube^n\to\F$ be non-zero. Suppose that $c_S\ne0$ in \eqref{eq:multilinear} and that $|S|=\deg_{\F}(f)$. Then $\supp(f)$ shatters $S^c=[n]\setminus S$.
\end{theorem}
Although its proof is elementary, this observation reveals a crucial algebraic viewpoint that considerably simplifies many subsequent arguments. In particular, the numerical uncertainty principle follows immediately.
\begin{corollary}[VC-dimension vs degree]\label{cor:main-ineq}
Let $\F$ be a field and let $f:\cube^n\to\F$ be non-zero.  Then
\begin{equation}\label{eq:VC-degree-field}
    \VC(f)+\deg_{\F}(f)\ge n.
\end{equation}
In particular, every non-zero Boolean function $f:\cube^n\to\cube$ satisfies
\begin{equation}\label{eq:VC-degree-bool}
    \VC(f)+\deg(f)\ge n,
    \qquad
    \VC(f)+\deg_{\Ftwo}(f)\ge n.
\end{equation}
\end{corollary}
\begin{remark}
This improves the trade-off obtained by combining two classical estimates.  If $\mathcal F\subseteq2^{[n]}$ and $\VC(\mathcal F)\le d$, then the Sauer--Shelah--Perles--Vapnik--Chervonenkis lemma \cite{Sauer1972,Shelah1972,VC1971} gives
\[
|\mathcal F|\le \sum_{i=0}^d \binom{n}{i}\le \left(\frac{en}{d}\right)^d.
\]
A standard Schwartz--Zippel type bound on the cube says that a real multilinear polynomial $g$ of degree $d^*$ satisfies
\[
    |\supp(g)|\ge 2^{n-d^*};
\]
see~\cite[Lemma 2.6]{NS1994degree}.  Applying these two estimates to $g=\mathbbm{1}_{\mathcal F}$ yields only
\[
d\log_2\left(\frac{en}{d}\right)+d^*\ge n.
\]
By contrast, \cref{cor:main-ineq} gives the sharp linear bound $d+d^*\ge n$. 
\end{remark}
\medskip
\noindent\textbf{Fourier support.}
For a Boolean function $f:\cube^n\to\cube$, let
\[
\Spec(f):=\{S\subseteq[n]:\hat f(S)\ne0\}
\]
be its real Fourier support. Bernasconi and Codenotti~\cite[Lemma 3]{BC1999} observed that
\[
\deg_{\Ftwo}(f)\le \log_2|\Spec(f)|.
\]
Since $|\supp(f)|\ge2^{\VC(f)}$, \cref{cor:main-ineq} recovers the classical uncertainty principle
\[
\log_2|\supp(f)|+\log_2|\Spec(f)|\ge n
\]
for non-zero Boolean functions; see also \cite[Exercise 3.15]{Ryanbook2014} and~\cite[Claim 3.3]{GT2013}. We prove a stronger shattering statement
\begin{equation}\label{eq:fourier-VC-intro}
    \VC(\Spec(f))\ge \deg_{\Ftwo}(f),
    \qquad
    \VC(\supp(f))+\VC(\Spec(f))\ge n.
\end{equation}
Thus, bounded VC-dimension of the support forces the Fourier support to shatter many coordinates.

\medskip
\noindent\textbf{Low-degree obstructions and Sauer-type bounds.}
The same theorem gives a short proof of the Sauer--Shelah lemma. If $\mathcal{F}\subseteq2^{[n]}$ has $\VC(\mathcal{F})\le d$, then no non-zero function of degree at most $n-d-1$ can be supported in $\mathcal{F}$. Hence, the restriction map from degree-at-most-$n-d-1$ functions to $\F^{\cube^n\setminus\mathcal F}$ is injective.  Taking dimensions gives
\[
    |\mathcal F|\le \sum_{i=0}^d \binom{n}{i}.
\]
This is Sauer--Shelah--Perles--Vapnik--Chervonenkis lemma~\cite{Sauer1972,Shelah1972,VC1971}, recovered from a low-degree obstruction.

We also prove a finite-product analogue.  Let $X=X_1\times\cdots\times X_n$ be a finite product set. For $T\subseteq[n]$, write $X_T:=\prod_{i\in T}X_i$, with the convention that $X_\emptyset$ is a one-point set. Let
\[
\pi_T:X\to X_T,\qquad\pi_T(x_1,\ldots,x_n)=(x_i)_{i\in T},
\]
be the coordinate projection. If $\mathcal A\subseteq X$, then $\pi_T(\mathcal A):=\{\pi_T(x):x\in\mathcal A\}\subseteq X_T$. For a nonempty set $\mathcal A\subseteq X$, define its coordinate projection dimension by
\[
\pdim(\mathcal A):=\max\left\{|T|:\pi_T(\mathcal A)=X_T\right\}.
\]
Thus $\pdim(\mathcal A)$ is the largest number of coordinates on which
$\mathcal A$ realizes every possible coordinate pattern. This is the full-projection, or coordinate-density, notion of Karpovsky--Milman~\cite{KM1978}; in Alon's terminology~\cite{Alon1983}, $\mathcal A$ is $I$-dense exactly when $\pi_I(\mathcal A)=X_I$. In the homogeneous case $X_1=\cdots=X_n=[q]$, it is the $q$-ary VC-dimension~\cite{MU2002}, also called the Karpovsky--Milman dimension in later
work on multivalued shattering \cite{FS2012}. 

We use the standard degree associated with the Efron--Stein decomposition on product spaces~\cite[Section~8.3]{Ryanbook2014}. Namely, write
\[
f=\sum_{S\subseteq[n]} f^{=S}
\]
for the Efron--Stein degree decomposition.  Equivalently, after choosing complements $\mathbb F^{X_i}=\langle 1\rangle\oplus W_i$ to the constants, the component $f^{=S}$ lies in $\bigotimes_{i\in S}W_i\otimes\bigotimes_{i\notin S}\langle 1\rangle$. We define
\[
\deg_\otimes(f):=\max\{|S|:f^{=S}\ne0\}.
\]
Equivalently, $\deg_\otimes(f)\le r$ means that $f$ lies in the linear span of
all functions depending on at most $r$ coordinates. With this degree, we prove
\begin{equation}
\pdim(\supp(f))+\deg_\otimes(f)\ge n
\end{equation}
for every nonzero function $f:X_1\times\cdots\times X_n\to\mathbb F$.
\medskip
\noindent\textbf{Polynomial-method consequences.}
Sziklai and Weiner \cite{sziklai2023covering} proved a sharp lower bound for polynomials that vanish on all high-weight vertices of the Boolean cube.  We recover this bound over arbitrary fields, with a weaker nonvanishing hypothesis.

\begin{corollary}[Sziklai--Weiner over arbitrary fields]\label{cor:SW-intro}
Let $\F$ be a field, let $0\le r\le n$, and let $P\in\F[X_1,\ldots,X_n]$. Denote $w_H(x)$ the Hamming weight simply as the number of $1$'s in $x$. Suppose that $P(x)=0$ for every $x\in\cube^n$ with $w_H(x)>r$ and that $P(x_0)\ne0$ for some $x_0\in\cube^n$ with $w_H(x_0)\le r$. Then $\deg(P)\ge n-r$.
\end{corollary}

Indeed, the restriction of $P$ to the cube is supported in the Hamming ball of radius $r$, which has VC-dimension $r$, and reduction modulo $x_i^2=x_i$ does not increase degree. This recovers the Sziklai--Weiner lower bound, whose original statement assumes nonvanishing on all vertices of weight at most $r$; the present formulation only requires one such vertex. 

The same viewpoint gives short proofs of a null-design shattering theorem. Null designs were introduced by Frankl and Pach~\cite{FP1983}; our conclusion strengthens the usual support-size consequence to a shattering statement.

\medskip
\noindent\textbf{Equality cases and organization.}
The inequalities in~\eqref{eq:VC-degree-bool} are sharp. If $C\subseteq\cube^n$ is a subcube of codimension $k$, then
\[
\VC(C)=n-k,\qquad\deg(\mathbbm{1}_C)=\deg_{\Ftwo}(\mathbbm{1}_C)=k.
\]
However, subcubes are not the only equality cases; non-subcube examples already occur for $n=4$.  We end with examples and computational data for $n\le4$, identifying supports up to the natural hypercube symmetries.

\Cref{sec:main-proof} proves \cref{thm:top-monomial}.  \Cref{sec:VC-applications} develops the Sauer--Shelah proof, the product-space theorem, the multivalued Sauer lemma, and the finite-Abelian-group version.  \Cref{sec:fourier} proves the Fourier support theorem and the polynomial-method applications, including \cref{cor:SW-intro}.  \Cref{sec:equality} records equality examples and computational data.

\section{Proof of the main theorem}\label{sec:main-proof}

We first prove the structural form, from which all field-specific inequalities follow.

\begin{proof}[Proof of \cref{thm:top-monomial}]
Let $D=\deg_{\F}(f)$, and let $S\subseteq[n]$ be such that $|S|=D$ and $c_S\ne0$ in the multilinear representation~\eqref{eq:multilinear}.  We show that every trace on $S^c$ is realized by a point of $\supp(f)$.

Fix an arbitrary $A\subseteq S^c$. Restrict the coordinates in $S^c$ by setting $x_i=1$ for $i\in A$ and $x_i=0$ for $i\in S^c\setminus A$. Let
\[
    g_A:\cube^S\to\F
\]
be the resulting function on the remaining coordinates. In the multilinear polynomial for $g_A$, the coefficient of the monomial $x^S$ is still $c_S$. Indeed, a monomial of $f$ can contribute to $x^S$ after this restriction only if it contains every variable in $S$. If it also contained any variable outside $S$, then its degree in $f$ would be larger than $|S|=D$, which is impossible by the choice of $D$.  Hence no other monomial changes the coefficient of $x^S$.

Thus $g_A$ is not the zero function.  Therefore there is some $B\subseteq S$ such that
\[
g_A(\one_B)=f(\one_{A\cup B})\ne0.
\]
The set $A\cup B\in\supp(f)$ has trace $A$ on $S^c$. Since this holds for every $A\subseteq S^c$, the support of $f$ shatters $S^c$.
\end{proof}

\begin{proof}[Proof of \cref{cor:main-ineq}]
By~\cref{thm:top-monomial}, $\supp(f)$ shatters $S^c$ for some $S$ of size $\deg_{\F}(f)$. Hence
\[
    \VC(f)\ge |S^c|=n-\deg_{\F}(f),
\]
which is equivalent to~\eqref{eq:VC-degree-field}.
\end{proof}

\section{Consequences for VC theory and product spaces}\label{sec:VC-applications}

\subsection{A linear-algebraic proof of Sauer--Shelah}

We first show that \cref{cor:main-ineq} recovers the Sauer--Shelah--Perles--Vapnik--Chervonenkis lemma.

\begin{theorem}[Sauer--Shelah--Perles--Vapnik--Chervonenkis]\label{thm:sauer}
Let $\mathcal A\subseteq2^{[n]}$ satisfy $\VC(\mathcal A)\le d$.  Then
\[
    |\mathcal A|\le \sum_{i=0}^d \binom{n}{i}.
\]
\end{theorem}

\begin{proof}
The case $d=n$ is trivial, so assume $d<n$.  Let $V^{\le n-d-1}$ be the vector space of all functions $g:\cube^n\to\F$ of degree at most $n-d-1$. Then
\[
    \dim V^{\le n-d-1}=\sum_{i=0}^{n-d-1}\binom{n}{i}.
\]
Consider the restriction map
\[
    \rho:V^{\le n-d-1}\longrightarrow \F^{\cube^n\setminus\mathcal A},
    \qquad
    \rho(g)=g|_{\cube^n\setminus\mathcal A}.
\]
We claim that $\rho$ is injective. If not, then some non-zero $g\in V^{\le n-d-1}$ vanishes on $\cube^n\setminus\mathcal A$, and hence $\supp(g)\subseteq\mathcal A$. Therefore,
\[
    \VC(g)=\VC(\supp(g))\le \VC(\mathcal A)\le d.
\]
On the other hand, \cref{cor:main-ineq} gives
\[
\deg_{\F}(g)\ge n-\VC(g)\ge n-d,
\]
contradicting $\deg_{\F}(g)\le n-d-1$. Thus $\rho$ is injective, and so
\[
|\cube^n\setminus\mathcal A|\ge \dim V_{\le n-d-1}=\sum_{i=0}^{n-d-1}\binom{n}{i}.
\]
Consequently,
\[
|\mathcal A|\le 2^n-\sum_{i=0}^{n-d-1}\binom{n}{i}=\sum_{i=0}^d \binom{n}{i},
\]
using the symmetry of binomial coefficients.
\end{proof}

\begin{remark}
This proof is close in spirit to algebraic and standard-monomial proofs of Sauer-type inequalities; see, for example, the order-shattering/standard-monomial approach of Anstee, R\'onyai, and Sali \cite{ARS2002} and later Groebner-basis treatments of shattering-extremal set systems \cite{Meszaros2020,MeszarosRonyai2011}.  The difference is that the key obstruction is phrased directly as a support-shattering statement: a non-zero low-degree function cannot be supported on a family of too small VC-dimension.
\end{remark}

\subsection{A product-space uncertainty principle}
We use the notation introduced above: $X=X_1\times\cdots\times X_n$, $\pdim$ is the coordinate projection dimension, and
\[
f=\sum_{S\subseteq[n]} f^{=S}
\]
is the Efron--Stein decomposition associated with fixed complements $\F^{X_i}=\langle1\rangle\oplus W_i$. Thus $\deg_\otimes(f)=\max\{|S|:f^{=S}\ne0\}$.
\begin{theorem}[Product-space uncertainty principle]\label{thm:product}
Let $f:X_1\times\cdots\times X_n\to\F$ be non-zero.  If $S\subseteq[n]$ satisfies $f^{=S}\ne0$ and $|S|=\deg_\otimes(f)$, then
\[
    \pi_{S^c}(\supp(f))=X_{S^c}.
\]
In particular,
\[
    \pdim(\supp(f))+\deg_\otimes(f)\ge n.
\]
\end{theorem}

\begin{proof}
Set $d=\deg_\otimes(f)$, and choose $S\subseteq[n]$ with $|S|=d$ and $f^{=S}\ne0$. We show that every point of $X_{S^c}$ is realized as the $S^c$-projection of a point in $\supp(f)$.

Fix $a\in X_{S^c}$. Let $f_a:X_S\to\F$ be the restriction of $f$ to the fiber over $a$, namely
\[
f_a(x_S)=f(x_S,a).
\]
Consider the induced Efron--Stein decomposition on $X_S$. The restriction of $f^{=S}$ to this fiber is still the same non-zero component in $\bigotimes_{i\in S}W_i$. No other component of $f$ contributes to this full $S$-component after the restriction. Indeed, if $U\subseteq S$ and $U\ne S$, then $f^{=U}$ has degree $|U|<|S|$ on $X_S$. If $U\not\subseteq S$, then fixing the coordinates outside $S$ turns $f^{=U}$ into a function whose components on $X_S$ are indexed by subsets of $U\cap S$, and $|U\cap S|<|S|$ because $|U|\le d=|S|$.

Therefore $f_a$ has a non-zero $S$-component, so $f_a$ is not identically zero. Hence there is some $x_S\in X_S$ such that $f(x_S,a)\ne0$.  Since this holds for every $a\in X_{S^c}$, we have
\[
    \pi_{S^c}(\supp(f))=X_{S^c}.
\]
Thus $\pdim(\supp(f))\ge |S^c|=n-d$, which is the desired inequality.
\end{proof}

\subsection{Multivalued Sauer lemma}

We next recover the full-projection multivalued Sauer lemma of Karpovsky--Milman, in Alon's non-homogeneous coordinate-density form~\cite{KM1978,Alon1983}. Let $X=X_1\times\cdots\times X_n$ and put $q_i=|X_i|$.

\begin{theorem}[Multivalued Sauer lemma]\label{thm:multivalued-sauer}
Let $\mathcal A\subseteq X_1\times\cdots\times X_n$ satisfy $\pdim(\mathcal A)\le d$. Then
\begin{equation}\label{eq:multivalued-sauer}
|\mathcal A|\le\sum_{\substack{S\subseteq[n]\\ |S|\ge n-d}}\prod_{i\in S}(q_i-1).
\end{equation}
In particular, if $q_1=\cdots=q_n=q$, then
\[
|\mathcal A|\le\sum_{i=0}^d \binom{n}{i}(q-1)^{n-i}.
\]
\end{theorem}

\begin{proof}
The case $d=n$ is trivial, so assume $d<n$. Let $L$ be the vector space of functions $g:X\to\F$ with $\deg_\otimes(g)\le n-d-1$. From the Efron--Stein decomposition,
\[
\dim L=\sum_{\substack{S\subseteq[n]\\ |S|\le n-d-1}}\prod_{i\in S}(q_i-1),
\]
since the $S$-component has dimension $\prod_{i\in S}(q_i-1)$.

Consider the restriction map
\[
    \rho:L\longrightarrow \F^{X\setminus\mathcal A},
    \qquad
    \rho(g)=g|_{X\setminus\mathcal A}.
\]
We claim that $\rho$ is injective. Otherwise, some non-zero $g\in L$ vanishes on $X\setminus\mathcal A$, so $\supp(g)\subseteq\mathcal A$. Hence
\[
    \pdim(\supp(g))\le \pdim(\mathcal A)\le d.
\]
By \cref{thm:product},
\[
\deg_\otimes(g)\ge n-\pdim(\supp(g))\ge n-d,
\]
contradicting $g\in L$. Thus $\rho$ is injective, and therefore
\[
|X\setminus\mathcal A|\ge\sum_{\substack{S\subseteq[n]\\ |S|\le n-d-1}}\prod_{i\in S}(q_i-1).
\]
Since
\[
|X|=\prod_{i=1}^n q_i=\prod_{i=1}^n(1+(q_i-1))=\sum_{S\subseteq[n]}\prod_{i\in S}(q_i-1),
\]
subtracting the lower bound for $|X\setminus\mathcal A|$ gives~\eqref{eq:multivalued-sauer}.
\end{proof}

\subsection{Finite Abelian groups}

Let $G=G_1\times\cdots\times G_n$ be a product of finite Abelian groups.  Every character $\chi$ of $G$ factors as
\[
    \chi=\chi_1\otimes\cdots\otimes\chi_n,
\]
where $\chi_i$ is a character of $G_i$. Define the Fourier degree of a character by
\[
\deg(\chi):=|\{i:\chi_i\text{ is non-trivial}\}|,
\]
and for $f:G\to\C$, define
\[
\deg_{\hat{G}}(f):=\max\{\deg(\chi):\hat{f}(\chi)\ne0\}.
\]
\begin{corollary}\label{cor:abelian}
Let $G=G_1\times\cdots\times G_n$ be a product of finite Abelian groups, and let $f:G\to\C$ be non-zero.  Then
\[
\pdim(\supp(f))+\deg_{\hat G}(f)\ge n.
\]
\end{corollary}

\begin{proof}
For each $i$, take $W_i\subseteq\C^{G_i}$ to be the span of the non-trivial characters of $G_i$. Then
\[
\C^{G_i}=\langle1\rangle\oplus W_i.
\]
With this choice, the Efron--Stein component $f^{=S}$ is precisely the sum of the Fourier terms whose non-trivial coordinate factors occur exactly on $S$. Hence $\deg_\otimes(f)=\deg_{\hat{G}}(f)$. The result follows from~\cref{thm:product}.
\end{proof}

\section{Fourier support and polynomial-method applications}\label{sec:fourier}

\subsection{A Fourier support theorem}

We now return to Boolean functions and answer~\cref{prob:central} in a stronger form. For $f:\cube^n\to\cube$, recall that
\[
    \Spec(f)=\{S\subseteq[n]:\hat{f}(S)\ne0\}
\]
is the real Fourier support, viewed as a family of subsets of $[n]$.

\begin{theorem}[Fourier--VC uncertainty principle]\label{thm:fourier-VC}
Let $f:\cube^n\to\cube$ be non-zero. Then
\[
    \VC(\Spec(f))\ge \deg_{\Ftwo}(f).
\]
Consequently,
\[
\VC(\supp(f))+\VC(\Spec(f))\ge n.
\]
In particular, if $\VC(f)=d$, then $\Spec(f)$ shatters a set of size at least $n-d$.
\end{theorem}

\begin{proof}
Let $r=\deg_{\Ftwo}(f)$, and choose $T\subseteq[n]$ with $|T|=r$ such that $x^T$ has non-zero coefficient in the algebraic normal form of $f$ over $\Ftwo$. Define
\[
    g:\cube^T\to\cube,
    \qquad
    g(y)=f(y,0_{[n]\setminus T}).
\]
The coefficient of $x^T$ in the algebraic normal form of $g$ is still $1$. Hence
\[
    \bigoplus_{y\in\cube^T}g(y)=1,
\]
so $g$ has odd Hamming weight.

For every $R\subseteq T$, the Fourier coefficient of $g$ is
\[
    \hat{g}(R)=2^{-|T|}\sum_{y\in\cube^T}g(y)(-1)^{\sum_{i\in R}y_i}.
\]
The numerator is a sum of an odd number of signs $\pm1$, and therefore cannot be zero.  Thus $\hat{g}(R)\ne0$ for every $R\subseteq T$.

On the other hand, substituting $x_i=0$ for $i\notin T$ into the Fourier expansion of $f$ gives
\[
\hat g(R)=\sum_{\substack{S\subseteq[n]\\ S\cap T=R}}\hat f(S).
\]
Since $\hat g(R)\ne0$, there is at least one $S\in\Spec(f)$ with $S\cap T=R$. Hence $\Spec(f)$ shatters $T$, and so
\[
    \VC(\Spec(f))\ge |T|=\deg_{\Ftwo}(f).
\]
Combining this with~\cref{cor:main-ineq} over $\Ftwo$ gives the second assertion.
\end{proof}

\begin{corollary}[Classical uncertainty principle]\label{cor:classical-uncertainty}
Let $f:\cube^n\to\cube$ be non-zero.  Then
\[
    \log_2|\supp(f)|+\log_2|\Spec(f)|\ge n.
\]
\end{corollary}

\begin{proof}
By~\cref{thm:fourier-VC}, $\Spec(f)$ shatters a set of size at least $n-\VC(\supp(f))$. Hence
\[
|\Spec(f)|\ge 2^{n-\VC(\supp(f))}.
\]
Since $|\supp(f)|\ge2^{\VC(\supp(f))}$, the result follows.
\end{proof}

\begin{remark}
The statement $\VC(\Spec(f))\ge\deg(f)$ with real degree in place of algebraic degree is false in general: the parity function has real degree $n$ but Fourier support of size one.
\end{remark}

\subsection{Sziklai--Weiner and null designs}

The original Sziklai--Weiner theorem assumes that the polynomial is nonzero on all low-weight vertices.  The support-shattering proof only uses that the restriction to the cube is nonzero and supported in the Hamming ball.  This is why the nonvanishing hypothesis in~\cref{cor:SW-intro} can be weakened to a single low-weight vertex.

\begin{proof}[Proof of \cref{cor:SW-intro}]
Let $f:\cube^n\to\F$ be the restriction of $P$ to the Boolean cube. The assumption $P(x_0)\ne0$ makes $f$ non-zero, while the vanishing assumption gives
\[
    \supp(f)\subseteq\{x\in\cube^n:\|x\|_1\le r\}.
\]
This Hamming ball has VC-dimension $r$.  Reducing $P$ modulo the relations $x_i^2=x_i$ gives the unique multilinear representative of $f$ over $\F$ and does not increase degree.  Therefore
\[
\deg(P)\ge\deg_{\F}(f)\ge n-\VC(f)\ge n-r,
\]
by~\cref{cor:main-ineq}.
\end{proof}

\begin{remark}
The corollary needs no assumption on the characteristic of $\F$. It also slightly weakens the usual nonvanishing hypothesis: it is enough that $P$ is non-zero at one vertex of weight at most $r$.  The bound is sharp over every field, as witnessed by
\[
P(X_1,\ldots,X_n)=\prod_{i=1}^{n-r}(1-X_i),
\]
which has degree $n-r$, vanishes on every vertex of weight larger than $r$, and is non-zero at the origin.
\end{remark}

We next record the null-design consequence. Null designs were introduced by Frankl and Pach~\cite{FP1983} in connection with lower bounds on the number of sets in signed designs. Let $X$ be a finite set and let $g:2^X\to\F$. We say that $g$ is a \emph{null $d$-design over $\F$} if $g$ is non-zero and
\[
    \sum_{A\subseteq F\subseteq X} g(F)=0
\]
for every $A\subseteq X$ with $|A|\le d$.

\begin{corollary}[Null designs]\label{cor:null-design}
If $g:2^X\to\F$ is a nonzero null $d$-design, then $\supp(g)$ shatters a set
of size at least $d+1$. In particular,
\[
    \VC(g)\ge d+1.
\]
\end{corollary}

\begin{proof}
Put $n=|X|$ and identify $2^X$ with the Boolean cube on the ground set $X$.
Define
\[
h(A):=(-1)^{|A|}g(X\setminus A),\qquad A\subseteq X.
\]
Then $h\ne0$. Moreover,
\[
\supp(h)=\{X\setminus F:F\in\supp(g)\}.
\]
Taking complements preserves shattering: a family $\mathcal A\subseteq 2^X$ shatters $S\subseteq X$ if and only if $\{X\setminus A:A\in\mathcal A\}$ shatters $S$. Hence $\VC(h)=\VC(g)$.

Let
\[
    h=\sum_{T\subseteq X}c_T\prod_{i\in T}x_i
\]
be the multilinear representation of $h$. By Möbius inversion,
\[
c_T=\sum_{A\subseteq T}(-1)^{|T|-|A|}h(A)=(-1)^{|T|}\sum_{A\subseteq T}g(X\setminus A).
\]
Equivalently,
\[
c_T=(-1)^{|T|}\sum_{X\setminus T\subseteq F\subseteq X}g(F).
\]
If $|T|\ge n-d$, then $|X\setminus T|\le d$, so the last sum is zero by the null $d$-design condition. Therefore $c_T=0$ for every $T$ with
$|T|\ge n-d$. Since $h\ne0$, this in particular forces $d<n$, and
\[
\deg_{\F}(h)\le n-d-1.
\]
Applying~\cref{cor:main-ineq} to $h$ gives
\[
\VC(g)=\VC(h)\ge n-\deg_{\F}(h)\ge d+1.
\]
Thus $\supp(g)$ shatters a set of size at least $d+1$.
\end{proof}

\begin{remark}
Since a family shattering a set of size $d+1$ has at least $2^{d+1}$ support points, \cref{cor:null-design} recovers the usual Frankl--Pach type lower bound on the support size of a non-zero null $d$-design, and gives the stronger conclusion that the support contains all traces on some $(d+1)$-set.
\end{remark}

\section{Equality cases and computations}\label{sec:equality}

A natural question is to determine when equality holds in the Boolean cases
\[
\VC(f)+\deg(f)=n,\qquad \VC(f)+\deg_{\Ftwo}(f)=n.
\]
Subcubes give the first family of examples. A subcube of codimension $k$ is a set of the form
\[
C=\{x\in\cube^n:x_{i_1}=a_1,\ldots,x_{i_k}=a_k\},
\]
where $i_1<\cdots<i_k$ and $a_1,\ldots,a_k\in\cube$. Its indicator satisfies
\[
\VC(C)=n-k,\qquad\deg(\mathbbm{1}_C)=\deg_{\Ftwo}(\one_C)=k,
\]
so equality holds in both inequalities.

Subcubes are not the only equality cases.  For instance, when $n=4$, the Boolean function
\[
    f(x_1,x_2,x_3,x_4)
    =x_1+x_2+x_3x_4-x_1x_2-x_1x_4-x_2x_3
\]
has $\VC(f)=\deg(f)=2$, and hence $\VC(f)+\deg(f)=4$.  Its support is
\[
\begin{split}
\supp(f)=\{&(0,0,1,1),(0,1,0,0),(0,1,0,1),(0,1,1,1),\\
&(1,0,0,0),(1,0,1,0),(1,0,1,1),(1,1,0,0)\},
\end{split}
\]
which is not a subcube.

We enumerated equality cases for $n\le4$.  Two supports are identified if one can be obtained from the other by a symmetry of the hypercube, namely by permuting coordinates and complementing arbitrary coordinates.  The orbit counts are shown in \cref{tab:equality}.

\begin{table}[h]
\centering
\begin{tabular}{rrrrrr}
\toprule
$n$ & total functions & real equality & real orbits & $\Ftwo$ equality & $\Ftwo$ orbits \\
\midrule
1 & 4 & 3 & 2 & 3 & 2 \\
2 & 16 & 9 & 3 & 11 & 4 \\
3 & 256 & 55 & 6 & 83 & 10 \\
4 & 65536 & 633 & 12 & 2491 & 38 \\
\bottomrule
\end{tabular}
\caption{Equality cases for the real-degree and $\Ftwo$-degree inequalities.  Orbits are taken under the full automorphism group of the hypercube.  The zero function is included in the total count but not in the equality counts.}
\label{tab:equality}
\end{table}

The orbit counts show that equality is already structurally richer than the subcube examples in dimension four.  To make the computation reproducible at the structural level, \cref{app:orbits} lists orbit representatives for $n=4$.  A full characterization appears to require more than the top-monomial shattering condition alone.

\begin{problem}\label{prob:equality}
Characterize the Boolean functions $f:\cube^n\to\cube$ for which equality holds in either
\[
\VC(f)+\deg(f)=n\qquad\text{or}\qquad\VC(f)+\deg_{\Ftwo}(f)=n.
\]
\end{problem}

The code accompanying the enumeration is available at
\[
    \texttt{https://github.com/FangYijia/deg-VC}.
\]

\appendix
\section{Orbit representatives in dimension four}\label{app:orbits}

This appendix lists representatives for the equality orbits in dimension four.  A support is written as a set of bitstrings $x_1x_2x_3x_4$.  The acting group is the full automorphism group of the four-dimensional cube, generated by coordinate permutations and coordinate complements.  The representatives were chosen canonically by lexicographic minimization over the orbit; the choice of representative is not meant to encode additional structure.

\begin{longtable}{@{}rrrr>{\ttfamily\footnotesize\raggedright\arraybackslash}p{0.58\textwidth}@{}}
\caption{Representatives for the 12 real-degree equality orbits in dimension four.}\label{tab:real-reps}\\
\toprule
No. & $|\supp(f)|$ & $\VC(f)$ & $\deg(f)$ & support representative \\
\midrule
\endfirsthead
\toprule
No. & $|\supp(f)|$ & $\VC(f)$ & $\deg(f)$ & support representative \\
\midrule
\endhead
1 & 1 & 0 & 4 & 0000 \\
2 & 2 & 1 & 3 & 0000, 1000 \\
3 & 2 & 1 & 3 & 0010, 1100 \\
4 & 4 & 2 & 2 & 0000, 0100, 1000, 1100 \\
5 & 4 & 1 & 3 & 0000, 0010, 1000, 1100 \\
6 & 4 & 1 & 3 & 0000, 0001, 0010, 1100 \\
7 & 4 & 1 & 3 & 0000, 0001, 1100, 1110 \\
8 & 4 & 2 & 2 & 0001, 0110, 1001, 1110 \\
9 & 8 & 3 & 1 & 0000, 0010, 0100, 0110, 1000, 1010, 1100, 1110 \\
10 & 8 & 2 & 2 & 0000, 0001, 0100, 0110, 1000, 1001, 1100, 1110 \\
11 & 8 & 2 & 2 & 0001, 0010, 0011, 0110, 1001, 1100, 1101, 1110 \\
12 & 16 & 4 & 0 & 0000, 0001, 0010, 0011, 0100, 0101, 0110, 0111, 1000, 1001, 1010, 1011, 1100, 1101, 1110, 1111 \\
\bottomrule
\end{longtable}

\begin{longtable}{@{}rrrr>{\ttfamily\footnotesize\raggedright\arraybackslash}p{0.58\textwidth}@{}}
\caption{Representatives for the 38 $\Ftwo$-degree equality orbits in dimension four.}\label{tab:f2-reps}\\
\toprule
No. & $|\supp(f)|$ & $\VC(f)$ & $\deg_{\Ftwo}(f)$ & support representative \\
\midrule
\endfirsthead
\toprule
No. & $|\supp(f)|$ & $\VC(f)$ & $\deg_{\Ftwo}(f)$ & support representative \\
\midrule
\endhead
1 & 1 & 0 & 4 & 0000 \\
2 & 2 & 1 & 3 & 0000, 1000 \\
3 & 2 & 1 & 3 & 0100, 1000 \\
4 & 2 & 1 & 3 & 0010, 1100 \\
5 & 2 & 1 & 3 & 0001, 1110 \\
6 & 4 & 2 & 2 & 0000, 0100, 1000, 1100 \\
7 & 4 & 1 & 3 & 0000, 0010, 0100, 1000 \\
8 & 4 & 1 & 3 & 0000, 0010, 1000, 1100 \\
9 & 4 & 2 & 2 & 0010, 0100, 1010, 1100 \\
10 & 4 & 2 & 2 & 0000, 0110, 1010, 1100 \\
11 & 4 & 1 & 3 & 0001, 0010, 0100, 1000 \\
12 & 4 & 1 & 3 & 0000, 0001, 0010, 1100 \\
13 & 4 & 1 & 3 & 0001, 0010, 1000, 1100 \\
14 & 4 & 1 & 3 & 0000, 0001, 1000, 1110 \\
15 & 4 & 1 & 3 & 0000, 0001, 1100, 1110 \\
16 & 4 & 2 & 2 & 0001, 0110, 1001, 1110 \\
17 & 4 & 2 & 2 & 0101, 0110, 1001, 1010 \\
18 & 4 & 2 & 2 & 0010, 0101, 1001, 1110 \\
19 & 6 & 2 & 2 & 0001, 0010, 0100, 0110, 1000, 1001 \\
20 & 6 & 2 & 2 & 0000, 0001, 0010, 0110, 1001, 1100 \\
21 & 6 & 2 & 2 & 0001, 0110, 1000, 1001, 1010, 1100 \\
22 & 6 & 2 & 2 & 0000, 0001, 0011, 0101, 1001, 1110 \\
23 & 8 & 3 & 1 & 0000, 0010, 0100, 0110, 1000, 1010, 1100, 1110 \\
24 & 8 & 2 & 2 & 0000, 0001, 0010, 0100, 1000, 1001, 1010, 1100 \\
25 & 8 & 2 & 2 & 0000, 0001, 0100, 0110, 1000, 1001, 1100, 1110 \\
26 & 8 & 2 & 2 & 0000, 0100, 0101, 0110, 1000, 1001, 1010, 1100 \\
27 & 8 & 2 & 2 & 0000, 0010, 0100, 0101, 1000, 1001, 1100, 1110 \\
28 & 8 & 2 & 2 & 0010, 0100, 0101, 0110, 1000, 1001, 1010, 1110 \\
29 & 8 & 3 & 1 & 0001, 0010, 0101, 0110, 1001, 1010, 1101, 1110 \\
30 & 8 & 2 & 2 & 0001, 0011, 0100, 0101, 1000, 1001, 1100, 1110 \\
31 & 8 & 2 & 2 & 0001, 0010, 0011, 0110, 1001, 1100, 1101, 1110 \\
32 & 8 & 3 & 1 & 0000, 0011, 0101, 0110, 1000, 1011, 1101, 1110 \\
33 & 8 & 3 & 1 & 0001, 0010, 0100, 0111, 1000, 1011, 1101, 1110 \\
34 & 10 & 2 & 2 & 0001, 0010, 0011, 0100, 0101, 0110, 1000, 1001, 1010, 1100 \\
35 & 10 & 2 & 2 & 0000, 0001, 0010, 0011, 0100, 0101, 1001, 1010, 1100, 1110 \\
36 & 10 & 2 & 2 & 0000, 0001, 0010, 0011, 0100, 0110, 1000, 1001, 1101, 1110 \\
37 & 10 & 2 & 2 & 0011, 0100, 0101, 0110, 1000, 1001, 1010, 1100, 1101, 1110 \\
38 & 16 & 4 & 0 & 0000, 0001, 0010, 0011, 0100, 0101, 0110, 0111, 1000, 1001, 1010, 1011, 1100, 1101, 1110, 1111 \\
\bottomrule
\end{longtable}

\section*{Acknowledgement}
The authors would like to thank the anonymous referees for their careful reading of the manuscript and for many valuable comments and suggestions, which significantly improved the quality of the paper. The authors are especially grateful to one anonymous referee whose helpful suggestions enabled the authors to simplify the proof.

This work was initiated during the 3$^{\rm rd}$ ECOPRO Student Research Program at the Institute for Basic Science (IBS) in the summer of 2025.  Chang and Fang are grateful to Professor Hong Liu for providing this research opportunity. The authors also thank Professor Hao Huang for helpful comments on an earlier draft of this manuscript.

The work of Fan Chang was supported by the National Natural Science Foundation of China (NSFC) under grant 124B2019 and by the Institute for Basic Science (IBS-R029-C4).

\bibliographystyle{abbrv}
\bibliography{reference}
\end{document}